\documentclass[12pt,a4paper]{amsart}

\usepackage{amsaddr}
\usepackage{fullpage}

\usepackage{scalerel,amssymb}
\usepackage{amsthm,aliascnt}
\usepackage{nicefrac}
\usepackage{amsmath}
\makeatletter
\renewcommand{\email}[2][]{%
  \ifx\emails\@empty\relax\else{\g@addto@macro\emails{,\space}}\fi%
  \@ifnotempty{#1}{\g@addto@macro\emails{\textrm{(#1)}\space}}%
  \g@addto@macro\emails{#2}%
}
\makeatother

\usepackage[utf8]{inputenc}
\usepackage[T1]{fontenc}
\usepackage[british]{babel}

\usepackage[
backend=bibtex,
doi=false,url=false,
maxbibnames=10 
]{biblatex}
\usepackage{xcolor}
\usepackage[breaklinks]{hyperref}
\usepackage[capitalize,nameinlink,noabbrev]{cleveref}
\usepackage{breakurl}

\definecolor{blackred}{RGB}{183, 24, 82}
\definecolor{lgreen}{rgb}{0.0, 0.48, 0.0}
\definecolor{lpurple}{rgb}{0.48, 0.0, 0.48}
\definecolor{bblue}{rgb}{0.2, 0.4, 0.8}
\hypersetup{linktocpage,
            colorlinks=true,
            linkcolor=lgreen,
            citecolor=lpurple,
            linktoc=true}
\makeatletter
\renewcommand{\tocsection}[3]{%
  \indentlabel{\@ifnotempty{#2}{\bfseries\ignorespaces#1 #2\quad}}\bfseries#3}
\renewcommand{\tocsubsection}[3]{%
  \indentlabel{\@ifnotempty{#2}{\ignorespaces#1 #2\quad}}#3}
\makeatother

\usepackage{floatrow}
\usepackage{float}
\usepackage{algorithm,algorithmic}
\Crefname{ALC@unique}{Line}{Lines}

\usepackage{graphicx}
\usepackage{epstopdf}
\usepackage{float}
\graphicspath{{./images/}}
\usepackage{caption, subcaption}

\usepackage{enumerate}
\usepackage{tabu}
\usepackage{multirow}

\usepackage[draft]{minted}

\makeatletter
\def\PYG@reset{\let\PYG@it=\relax \let\PYG@bf=\relax%
    \let\PYG@ul=\relax \let\PYG@tc=\relax%
    \let\PYG@bc=\relax \let\PYG@ff=\relax}
\def\PYG@tok#1{\csname PYG@tok@#1\endcsname}
\def\PYG@toks#1+{\ifx\relax#1\empty\else%
    \PYG@tok{#1}\expandafter\PYG@toks\fi}
\def\PYG@do#1{\PYG@bc{\PYG@tc{\PYG@ul{%
    \PYG@it{\PYG@bf{\PYG@ff{#1}}}}}}}
\def\PYG#1#2{\PYG@reset\PYG@toks#1+\relax+\PYG@do{#2}}

\expandafter\def\csname PYG@tok@w\endcsname{\def\PYG@tc##1{\textcolor[rgb]{0.73,0.73,0.73}{##1}}}
\expandafter\def\csname PYG@tok@c\endcsname{\let\PYG@it=\textit\def\PYG@tc##1{\textcolor[rgb]{0.25,0.50,0.50}{##1}}}
\expandafter\def\csname PYG@tok@cp\endcsname{\def\PYG@tc##1{\textcolor[rgb]{0.74,0.48,0.00}{##1}}}
\expandafter\def\csname PYG@tok@k\endcsname{\let\PYG@bf=\textbf\def\PYG@tc##1{\textcolor[rgb]{0.00,0.50,0.00}{##1}}}
\expandafter\def\csname PYG@tok@kp\endcsname{\def\PYG@tc##1{\textcolor[rgb]{0.00,0.50,0.00}{##1}}}
\expandafter\def\csname PYG@tok@kt\endcsname{\def\PYG@tc##1{\textcolor[rgb]{0.69,0.00,0.25}{##1}}}
\expandafter\def\csname PYG@tok@o\endcsname{\def\PYG@tc##1{\textcolor[rgb]{0.40,0.40,0.40}{##1}}}
\expandafter\def\csname PYG@tok@ow\endcsname{\let\PYG@bf=\textbf\def\PYG@tc##1{\textcolor[rgb]{0.67,0.13,1.00}{##1}}}
\expandafter\def\csname PYG@tok@nb\endcsname{\def\PYG@tc##1{\textcolor[rgb]{0.00,0.50,0.00}{##1}}}
\expandafter\def\csname PYG@tok@nf\endcsname{\def\PYG@tc##1{\textcolor[rgb]{0.00,0.00,1.00}{##1}}}
\expandafter\def\csname PYG@tok@nc\endcsname{\let\PYG@bf=\textbf\def\PYG@tc##1{\textcolor[rgb]{0.00,0.00,1.00}{##1}}}
\expandafter\def\csname PYG@tok@nn\endcsname{\let\PYG@bf=\textbf\def\PYG@tc##1{\textcolor[rgb]{0.00,0.00,1.00}{##1}}}
\expandafter\def\csname PYG@tok@ne\endcsname{\let\PYG@bf=\textbf\def\PYG@tc##1{\textcolor[rgb]{0.82,0.25,0.23}{##1}}}
\expandafter\def\csname PYG@tok@nv\endcsname{\def\PYG@tc##1{\textcolor[rgb]{0.10,0.09,0.49}{##1}}}
\expandafter\def\csname PYG@tok@no\endcsname{\def\PYG@tc##1{\textcolor[rgb]{0.53,0.00,0.00}{##1}}}
\expandafter\def\csname PYG@tok@nl\endcsname{\def\PYG@tc##1{\textcolor[rgb]{0.63,0.63,0.00}{##1}}}
\expandafter\def\csname PYG@tok@ni\endcsname{\let\PYG@bf=\textbf\def\PYG@tc##1{\textcolor[rgb]{0.60,0.60,0.60}{##1}}}
\expandafter\def\csname PYG@tok@na\endcsname{\def\PYG@tc##1{\textcolor[rgb]{0.49,0.56,0.16}{##1}}}
\expandafter\def\csname PYG@tok@nt\endcsname{\let\PYG@bf=\textbf\def\PYG@tc##1{\textcolor[rgb]{0.00,0.50,0.00}{##1}}}
\expandafter\def\csname PYG@tok@nd\endcsname{\def\PYG@tc##1{\textcolor[rgb]{0.67,0.13,1.00}{##1}}}
\expandafter\def\csname PYG@tok@s\endcsname{\def\PYG@tc##1{\textcolor[rgb]{0.73,0.13,0.13}{##1}}}
\expandafter\def\csname PYG@tok@sd\endcsname{\let\PYG@it=\textit\def\PYG@tc##1{\textcolor[rgb]{0.73,0.13,0.13}{##1}}}
\expandafter\def\csname PYG@tok@si\endcsname{\let\PYG@bf=\textbf\def\PYG@tc##1{\textcolor[rgb]{0.73,0.40,0.53}{##1}}}
\expandafter\def\csname PYG@tok@se\endcsname{\let\PYG@bf=\textbf\def\PYG@tc##1{\textcolor[rgb]{0.73,0.40,0.13}{##1}}}
\expandafter\def\csname PYG@tok@sr\endcsname{\def\PYG@tc##1{\textcolor[rgb]{0.73,0.40,0.53}{##1}}}
\expandafter\def\csname PYG@tok@ss\endcsname{\def\PYG@tc##1{\textcolor[rgb]{0.10,0.09,0.49}{##1}}}
\expandafter\def\csname PYG@tok@sx\endcsname{\def\PYG@tc##1{\textcolor[rgb]{0.00,0.50,0.00}{##1}}}
\expandafter\def\csname PYG@tok@m\endcsname{\def\PYG@tc##1{\textcolor[rgb]{0.40,0.40,0.40}{##1}}}
\expandafter\def\csname PYG@tok@gh\endcsname{\let\PYG@bf=\textbf\def\PYG@tc##1{\textcolor[rgb]{0.00,0.00,0.50}{##1}}}
\expandafter\def\csname PYG@tok@gu\endcsname{\let\PYG@bf=\textbf\def\PYG@tc##1{\textcolor[rgb]{0.50,0.00,0.50}{##1}}}
\expandafter\def\csname PYG@tok@gd\endcsname{\def\PYG@tc##1{\textcolor[rgb]{0.63,0.00,0.00}{##1}}}
\expandafter\def\csname PYG@tok@gi\endcsname{\def\PYG@tc##1{\textcolor[rgb]{0.00,0.63,0.00}{##1}}}
\expandafter\def\csname PYG@tok@gr\endcsname{\def\PYG@tc##1{\textcolor[rgb]{1.00,0.00,0.00}{##1}}}
\expandafter\def\csname PYG@tok@ge\endcsname{\let\PYG@it=\textit}
\expandafter\def\csname PYG@tok@gs\endcsname{\let\PYG@bf=\textbf}
\expandafter\def\csname PYG@tok@gp\endcsname{\let\PYG@bf=\textbf\def\PYG@tc##1{\textcolor[rgb]{0.00,0.00,0.50}{##1}}}
\expandafter\def\csname PYG@tok@go\endcsname{\def\PYG@tc##1{\textcolor[rgb]{0.53,0.53,0.53}{##1}}}
\expandafter\def\csname PYG@tok@gt\endcsname{\def\PYG@tc##1{\textcolor[rgb]{0.00,0.27,0.87}{##1}}}
\expandafter\def\csname PYG@tok@err\endcsname{\def\PYG@bc##1{\setlength{\fboxsep}{0pt}\fcolorbox[rgb]{1.00,0.00,0.00}{1,1,1}{\strut ##1}}}
\expandafter\def\csname PYG@tok@kc\endcsname{\let\PYG@bf=\textbf\def\PYG@tc##1{\textcolor[rgb]{0.00,0.50,0.00}{##1}}}
\expandafter\def\csname PYG@tok@kd\endcsname{\let\PYG@bf=\textbf\def\PYG@tc##1{\textcolor[rgb]{0.00,0.50,0.00}{##1}}}
\expandafter\def\csname PYG@tok@kn\endcsname{\let\PYG@bf=\textbf\def\PYG@tc##1{\textcolor[rgb]{0.00,0.50,0.00}{##1}}}
\expandafter\def\csname PYG@tok@kr\endcsname{\let\PYG@bf=\textbf\def\PYG@tc##1{\textcolor[rgb]{0.00,0.50,0.00}{##1}}}
\expandafter\def\csname PYG@tok@bp\endcsname{\def\PYG@tc##1{\textcolor[rgb]{0.00,0.50,0.00}{##1}}}
\expandafter\def\csname PYG@tok@fm\endcsname{\def\PYG@tc##1{\textcolor[rgb]{0.00,0.00,1.00}{##1}}}
\expandafter\def\csname PYG@tok@vc\endcsname{\def\PYG@tc##1{\textcolor[rgb]{0.10,0.09,0.49}{##1}}}
\expandafter\def\csname PYG@tok@vg\endcsname{\def\PYG@tc##1{\textcolor[rgb]{0.10,0.09,0.49}{##1}}}
\expandafter\def\csname PYG@tok@vi\endcsname{\def\PYG@tc##1{\textcolor[rgb]{0.10,0.09,0.49}{##1}}}
\expandafter\def\csname PYG@tok@vm\endcsname{\def\PYG@tc##1{\textcolor[rgb]{0.10,0.09,0.49}{##1}}}
\expandafter\def\csname PYG@tok@sa\endcsname{\def\PYG@tc##1{\textcolor[rgb]{0.73,0.13,0.13}{##1}}}
\expandafter\def\csname PYG@tok@sb\endcsname{\def\PYG@tc##1{\textcolor[rgb]{0.73,0.13,0.13}{##1}}}
\expandafter\def\csname PYG@tok@sc\endcsname{\def\PYG@tc##1{\textcolor[rgb]{0.73,0.13,0.13}{##1}}}
\expandafter\def\csname PYG@tok@dl\endcsname{\def\PYG@tc##1{\textcolor[rgb]{0.73,0.13,0.13}{##1}}}
\expandafter\def\csname PYG@tok@s2\endcsname{\def\PYG@tc##1{\textcolor[rgb]{0.73,0.13,0.13}{##1}}}
\expandafter\def\csname PYG@tok@sh\endcsname{\def\PYG@tc##1{\textcolor[rgb]{0.73,0.13,0.13}{##1}}}
\expandafter\def\csname PYG@tok@s1\endcsname{\def\PYG@tc##1{\textcolor[rgb]{0.73,0.13,0.13}{##1}}}
\expandafter\def\csname PYG@tok@mb\endcsname{\def\PYG@tc##1{\textcolor[rgb]{0.40,0.40,0.40}{##1}}}
\expandafter\def\csname PYG@tok@mf\endcsname{\def\PYG@tc##1{\textcolor[rgb]{0.40,0.40,0.40}{##1}}}
\expandafter\def\csname PYG@tok@mh\endcsname{\def\PYG@tc##1{\textcolor[rgb]{0.40,0.40,0.40}{##1}}}
\expandafter\def\csname PYG@tok@mi\endcsname{\def\PYG@tc##1{\textcolor[rgb]{0.40,0.40,0.40}{##1}}}
\expandafter\def\csname PYG@tok@il\endcsname{\def\PYG@tc##1{\textcolor[rgb]{0.40,0.40,0.40}{##1}}}
\expandafter\def\csname PYG@tok@mo\endcsname{\def\PYG@tc##1{\textcolor[rgb]{0.40,0.40,0.40}{##1}}}
\expandafter\def\csname PYG@tok@ch\endcsname{\let\PYG@it=\textit\def\PYG@tc##1{\textcolor[rgb]{0.25,0.50,0.50}{##1}}}
\expandafter\def\csname PYG@tok@cm\endcsname{\let\PYG@it=\textit\def\PYG@tc##1{\textcolor[rgb]{0.25,0.50,0.50}{##1}}}
\expandafter\def\csname PYG@tok@cpf\endcsname{\let\PYG@it=\textit\def\PYG@tc##1{\textcolor[rgb]{0.25,0.50,0.50}{##1}}}
\expandafter\def\csname PYG@tok@c1\endcsname{\let\PYG@it=\textit\def\PYG@tc##1{\textcolor[rgb]{0.25,0.50,0.50}{##1}}}
\expandafter\def\csname PYG@tok@cs\endcsname{\let\PYG@it=\textit\def\PYG@tc##1{\textcolor[rgb]{0.25,0.50,0.50}{##1}}}

\def\PYGZgt{\char`\>}

\def\PYGZhy{\char`\-}

\makeatother

\makeatletter
\def\PYGdefault@reset{\let\PYGdefault@it=\relax \let\PYGdefault@bf=\relax%
    \let\PYGdefault@ul=\relax \let\PYGdefault@tc=\relax%
    \let\PYGdefault@bc=\relax \let\PYGdefault@ff=\relax}
\def\PYGdefault@tok#1{\csname PYGdefault@tok@#1\endcsname}
\def\PYGdefault@toks#1+{\ifx\relax#1\empty\else%
    \PYGdefault@tok{#1}\expandafter\PYGdefault@toks\fi}
\def\PYGdefault@do#1{\PYGdefault@bc{\PYGdefault@tc{\PYGdefault@ul{%
    \PYGdefault@it{\PYGdefault@bf{\PYGdefault@ff{#1}}}}}}}
\def\PYGdefault#1#2{\PYGdefault@reset\PYGdefault@toks#1+\relax+\PYGdefault@do{#2}}

\expandafter\def\csname PYGdefault@tok@w\endcsname{\def\PYGdefault@tc##1{\textcolor[rgb]{0.73,0.73,0.73}{##1}}}
\expandafter\def\csname PYGdefault@tok@c\endcsname{\let\PYGdefault@it=\textit\def\PYGdefault@tc##1{\textcolor[rgb]{0.25,0.50,0.50}{##1}}}
\expandafter\def\csname PYGdefault@tok@cp\endcsname{\def\PYGdefault@tc##1{\textcolor[rgb]{0.74,0.48,0.00}{##1}}}
\expandafter\def\csname PYGdefault@tok@k\endcsname{\let\PYGdefault@bf=\textbf\def\PYGdefault@tc##1{\textcolor[rgb]{0.00,0.50,0.00}{##1}}}
\expandafter\def\csname PYGdefault@tok@kp\endcsname{\def\PYGdefault@tc##1{\textcolor[rgb]{0.00,0.50,0.00}{##1}}}
\expandafter\def\csname PYGdefault@tok@kt\endcsname{\def\PYGdefault@tc##1{\textcolor[rgb]{0.69,0.00,0.25}{##1}}}
\expandafter\def\csname PYGdefault@tok@o\endcsname{\def\PYGdefault@tc##1{\textcolor[rgb]{0.40,0.40,0.40}{##1}}}
\expandafter\def\csname PYGdefault@tok@ow\endcsname{\let\PYGdefault@bf=\textbf\def\PYGdefault@tc##1{\textcolor[rgb]{0.67,0.13,1.00}{##1}}}
\expandafter\def\csname PYGdefault@tok@nb\endcsname{\def\PYGdefault@tc##1{\textcolor[rgb]{0.00,0.50,0.00}{##1}}}
\expandafter\def\csname PYGdefault@tok@nf\endcsname{\def\PYGdefault@tc##1{\textcolor[rgb]{0.00,0.00,1.00}{##1}}}
\expandafter\def\csname PYGdefault@tok@nc\endcsname{\let\PYGdefault@bf=\textbf\def\PYGdefault@tc##1{\textcolor[rgb]{0.00,0.00,1.00}{##1}}}
\expandafter\def\csname PYGdefault@tok@nn\endcsname{\let\PYGdefault@bf=\textbf\def\PYGdefault@tc##1{\textcolor[rgb]{0.00,0.00,1.00}{##1}}}
\expandafter\def\csname PYGdefault@tok@ne\endcsname{\let\PYGdefault@bf=\textbf\def\PYGdefault@tc##1{\textcolor[rgb]{0.82,0.25,0.23}{##1}}}
\expandafter\def\csname PYGdefault@tok@nv\endcsname{\def\PYGdefault@tc##1{\textcolor[rgb]{0.10,0.09,0.49}{##1}}}
\expandafter\def\csname PYGdefault@tok@no\endcsname{\def\PYGdefault@tc##1{\textcolor[rgb]{0.53,0.00,0.00}{##1}}}
\expandafter\def\csname PYGdefault@tok@nl\endcsname{\def\PYGdefault@tc##1{\textcolor[rgb]{0.63,0.63,0.00}{##1}}}
\expandafter\def\csname PYGdefault@tok@ni\endcsname{\let\PYGdefault@bf=\textbf\def\PYGdefault@tc##1{\textcolor[rgb]{0.60,0.60,0.60}{##1}}}
\expandafter\def\csname PYGdefault@tok@na\endcsname{\def\PYGdefault@tc##1{\textcolor[rgb]{0.49,0.56,0.16}{##1}}}
\expandafter\def\csname PYGdefault@tok@nt\endcsname{\let\PYGdefault@bf=\textbf\def\PYGdefault@tc##1{\textcolor[rgb]{0.00,0.50,0.00}{##1}}}
\expandafter\def\csname PYGdefault@tok@nd\endcsname{\def\PYGdefault@tc##1{\textcolor[rgb]{0.67,0.13,1.00}{##1}}}
\expandafter\def\csname PYGdefault@tok@s\endcsname{\def\PYGdefault@tc##1{\textcolor[rgb]{0.73,0.13,0.13}{##1}}}
\expandafter\def\csname PYGdefault@tok@sd\endcsname{\let\PYGdefault@it=\textit\def\PYGdefault@tc##1{\textcolor[rgb]{0.73,0.13,0.13}{##1}}}
\expandafter\def\csname PYGdefault@tok@si\endcsname{\let\PYGdefault@bf=\textbf\def\PYGdefault@tc##1{\textcolor[rgb]{0.73,0.40,0.53}{##1}}}
\expandafter\def\csname PYGdefault@tok@se\endcsname{\let\PYGdefault@bf=\textbf\def\PYGdefault@tc##1{\textcolor[rgb]{0.73,0.40,0.13}{##1}}}
\expandafter\def\csname PYGdefault@tok@sr\endcsname{\def\PYGdefault@tc##1{\textcolor[rgb]{0.73,0.40,0.53}{##1}}}
\expandafter\def\csname PYGdefault@tok@ss\endcsname{\def\PYGdefault@tc##1{\textcolor[rgb]{0.10,0.09,0.49}{##1}}}
\expandafter\def\csname PYGdefault@tok@sx\endcsname{\def\PYGdefault@tc##1{\textcolor[rgb]{0.00,0.50,0.00}{##1}}}
\expandafter\def\csname PYGdefault@tok@m\endcsname{\def\PYGdefault@tc##1{\textcolor[rgb]{0.40,0.40,0.40}{##1}}}
\expandafter\def\csname PYGdefault@tok@gh\endcsname{\let\PYGdefault@bf=\textbf\def\PYGdefault@tc##1{\textcolor[rgb]{0.00,0.00,0.50}{##1}}}
\expandafter\def\csname PYGdefault@tok@gu\endcsname{\let\PYGdefault@bf=\textbf\def\PYGdefault@tc##1{\textcolor[rgb]{0.50,0.00,0.50}{##1}}}
\expandafter\def\csname PYGdefault@tok@gd\endcsname{\def\PYGdefault@tc##1{\textcolor[rgb]{0.63,0.00,0.00}{##1}}}
\expandafter\def\csname PYGdefault@tok@gi\endcsname{\def\PYGdefault@tc##1{\textcolor[rgb]{0.00,0.63,0.00}{##1}}}
\expandafter\def\csname PYGdefault@tok@gr\endcsname{\def\PYGdefault@tc##1{\textcolor[rgb]{1.00,0.00,0.00}{##1}}}
\expandafter\def\csname PYGdefault@tok@ge\endcsname{\let\PYGdefault@it=\textit}
\expandafter\def\csname PYGdefault@tok@gs\endcsname{\let\PYGdefault@bf=\textbf}
\expandafter\def\csname PYGdefault@tok@gp\endcsname{\let\PYGdefault@bf=\textbf\def\PYGdefault@tc##1{\textcolor[rgb]{0.00,0.00,0.50}{##1}}}
\expandafter\def\csname PYGdefault@tok@go\endcsname{\def\PYGdefault@tc##1{\textcolor[rgb]{0.53,0.53,0.53}{##1}}}
\expandafter\def\csname PYGdefault@tok@gt\endcsname{\def\PYGdefault@tc##1{\textcolor[rgb]{0.00,0.27,0.87}{##1}}}
\expandafter\def\csname PYGdefault@tok@err\endcsname{\def\PYGdefault@bc##1{\setlength{\fboxsep}{0pt}\fcolorbox[rgb]{1.00,0.00,0.00}{1,1,1}{\strut ##1}}}
\expandafter\def\csname PYGdefault@tok@kc\endcsname{\let\PYGdefault@bf=\textbf\def\PYGdefault@tc##1{\textcolor[rgb]{0.00,0.50,0.00}{##1}}}
\expandafter\def\csname PYGdefault@tok@kd\endcsname{\let\PYGdefault@bf=\textbf\def\PYGdefault@tc##1{\textcolor[rgb]{0.00,0.50,0.00}{##1}}}
\expandafter\def\csname PYGdefault@tok@kn\endcsname{\let\PYGdefault@bf=\textbf\def\PYGdefault@tc##1{\textcolor[rgb]{0.00,0.50,0.00}{##1}}}
\expandafter\def\csname PYGdefault@tok@kr\endcsname{\let\PYGdefault@bf=\textbf\def\PYGdefault@tc##1{\textcolor[rgb]{0.00,0.50,0.00}{##1}}}
\expandafter\def\csname PYGdefault@tok@bp\endcsname{\def\PYGdefault@tc##1{\textcolor[rgb]{0.00,0.50,0.00}{##1}}}
\expandafter\def\csname PYGdefault@tok@fm\endcsname{\def\PYGdefault@tc##1{\textcolor[rgb]{0.00,0.00,1.00}{##1}}}
\expandafter\def\csname PYGdefault@tok@vc\endcsname{\def\PYGdefault@tc##1{\textcolor[rgb]{0.10,0.09,0.49}{##1}}}
\expandafter\def\csname PYGdefault@tok@vg\endcsname{\def\PYGdefault@tc##1{\textcolor[rgb]{0.10,0.09,0.49}{##1}}}
\expandafter\def\csname PYGdefault@tok@vi\endcsname{\def\PYGdefault@tc##1{\textcolor[rgb]{0.10,0.09,0.49}{##1}}}
\expandafter\def\csname PYGdefault@tok@vm\endcsname{\def\PYGdefault@tc##1{\textcolor[rgb]{0.10,0.09,0.49}{##1}}}
\expandafter\def\csname PYGdefault@tok@sa\endcsname{\def\PYGdefault@tc##1{\textcolor[rgb]{0.73,0.13,0.13}{##1}}}
\expandafter\def\csname PYGdefault@tok@sb\endcsname{\def\PYGdefault@tc##1{\textcolor[rgb]{0.73,0.13,0.13}{##1}}}
\expandafter\def\csname PYGdefault@tok@sc\endcsname{\def\PYGdefault@tc##1{\textcolor[rgb]{0.73,0.13,0.13}{##1}}}
\expandafter\def\csname PYGdefault@tok@dl\endcsname{\def\PYGdefault@tc##1{\textcolor[rgb]{0.73,0.13,0.13}{##1}}}
\expandafter\def\csname PYGdefault@tok@s2\endcsname{\def\PYGdefault@tc##1{\textcolor[rgb]{0.73,0.13,0.13}{##1}}}
\expandafter\def\csname PYGdefault@tok@sh\endcsname{\def\PYGdefault@tc##1{\textcolor[rgb]{0.73,0.13,0.13}{##1}}}
\expandafter\def\csname PYGdefault@tok@s1\endcsname{\def\PYGdefault@tc##1{\textcolor[rgb]{0.73,0.13,0.13}{##1}}}
\expandafter\def\csname PYGdefault@tok@mb\endcsname{\def\PYGdefault@tc##1{\textcolor[rgb]{0.40,0.40,0.40}{##1}}}
\expandafter\def\csname PYGdefault@tok@mf\endcsname{\def\PYGdefault@tc##1{\textcolor[rgb]{0.40,0.40,0.40}{##1}}}
\expandafter\def\csname PYGdefault@tok@mh\endcsname{\def\PYGdefault@tc##1{\textcolor[rgb]{0.40,0.40,0.40}{##1}}}
\expandafter\def\csname PYGdefault@tok@mi\endcsname{\def\PYGdefault@tc##1{\textcolor[rgb]{0.40,0.40,0.40}{##1}}}
\expandafter\def\csname PYGdefault@tok@il\endcsname{\def\PYGdefault@tc##1{\textcolor[rgb]{0.40,0.40,0.40}{##1}}}
\expandafter\def\csname PYGdefault@tok@mo\endcsname{\def\PYGdefault@tc##1{\textcolor[rgb]{0.40,0.40,0.40}{##1}}}
\expandafter\def\csname PYGdefault@tok@ch\endcsname{\let\PYGdefault@it=\textit\def\PYGdefault@tc##1{\textcolor[rgb]{0.25,0.50,0.50}{##1}}}
\expandafter\def\csname PYGdefault@tok@cm\endcsname{\let\PYGdefault@it=\textit\def\PYGdefault@tc##1{\textcolor[rgb]{0.25,0.50,0.50}{##1}}}
\expandafter\def\csname PYGdefault@tok@cpf\endcsname{\let\PYGdefault@it=\textit\def\PYGdefault@tc##1{\textcolor[rgb]{0.25,0.50,0.50}{##1}}}
\expandafter\def\csname PYGdefault@tok@c1\endcsname{\let\PYGdefault@it=\textit\def\PYGdefault@tc##1{\textcolor[rgb]{0.25,0.50,0.50}{##1}}}
\expandafter\def\csname PYGdefault@tok@cs\endcsname{\let\PYGdefault@it=\textit\def\PYGdefault@tc##1{\textcolor[rgb]{0.25,0.50,0.50}{##1}}}

\makeatother

\usepackage{etoolbox}
\BeforeBeginEnvironment{minted}{\medskip}
\AfterEndEnvironment{minted}{\medskip}

\usepackage{catoptions}
\makeatletter

\def\Autoref#1{%
  \begingroup
  \edef\reserved@a{\cpttrimspaces{#1}}%
  \ifcsndefTF{r@#1}{%
    \xaftercsname{\expandafter\testreftype\@fourthoffive}
      {r@\reserved@a}.\\{#1}%
  }{%
    \ref{#1}%
  }%
  \endgroup
}
\def\testreftype#1.#2\\#3{%
  \ifcsndefTF{#1autorefname}{%
    \def\reserved@a##1##2\@nil{%
      \uppercase{\def\ref@name{##1}}%
      \csn@edef{#1autorefname}{\ref@name##2}%
      \autoref{#3}%
    }%
    \reserved@a#1\@nil
  }{%
    \autoref{#3}%
  }%
}
\makeatother

\newcommand{\seq}[1]{\left(#1\right)}
\newcommand{\idx}[1]{\mbox{\underline{\sf #1}}}

\DeclareMathOperator*{\Seq}{\mbox{\sc Seq}}

\newcommand{\CS}[1]{\mathcal{#1}}

\let\emptyset\varnothing

\definecolor{lgreen}{rgb}{0.0, 0.48, 0.0}
\definecolor{lpurple}{rgb}{0.48, 0.0, 0.48}
\definecolor{bblue}{rgb}{0.2, 0.4, 0.8}
\hypersetup{linktocpage,
            colorlinks=true,
            linkcolor=lgreen,
            citecolor=lpurple,
            linktoc=true}

\definecolor{bblue}{rgb}{0.2, 0.4, 0.8}
\definecolor{bgreen}{rgb}{0.2, 0.6, 0.4}
\definecolor{bred}{rgb}{0.8, 0.4, 0.2}
\definecolor{bviolet}{rgb}{0.7, 0.2, 0.7}
\definecolor{blackred}{rgb}{0.6, 0.3, 0.3}
\definecolor{blackblue}{rgb}{0.3, 0.3, 0.6}

\usepackage{tikz}
\usetikzlibrary{fit,arrows,trees,shapes,shapes.geometric,calc,matrix}
\tikzset{
  treenode/.style = {align=center, inner sep=0pt, text centered,
    font=\sffamily},
  arn_nn/.style = {treenode, circle, bblue, draw=bblue,
    fill=bblue!10,
    minimum width=0.5em, minimum height=0.5em
},
  arn_n/.style = {treenode, circle, bblue, draw=bblue,
    text width=1.4em, very thick,
    fill=bblue!10},
  arn_g/.style = {treenode, circle, bgreen, draw=bgreen,
    minimum width=0.5em, minimum height=0.5em, text width=1.2em,
    fill=bblue!10},
  arn_r/.style = {treenode, circle, bred, draw=bred,
    minimum width=0.5em, minimum height=0.5em, text width=1.4em,
    fill=bviolet!10},
  arn_rp/.style = {treenode, circle, bred, draw=bred,
    minimum width=0.5em, minimum height=0.5em,
    fill=bviolet!10},
  arn_x/.style = {treenode, triangle, draw=black,
    minimum width=0.5em, minimum height=0.5em},
  triangle/.style = {treenode, bred, draw=bred, fill=bred!20, regular polygon, regular polygon
    sides=3, very thick, text width=1.5em },
  triangle_b/.style = {treenode, bblue, draw=bblue,
    fill=bblue!20, regular polygon, regular polygon
    sides=3, very thick, text width=1.5em },
  triangle_g/.style = {treenode, bgreen, draw=bgreen,
    fill=bgreen!20, regular polygon, regular polygon
    sides=3, very thick, text width=1.5em },
  triangle_v/.style = {treenode, bviolet, draw=bviolet,
    fill=bviolet!20, regular polygon, regular polygon
    sides=3, very thick, text width=1.5em },
  triangle_h/.style = {treenode, bblue, draw=bblue,
    fill=gray!20, regular polygon, regular polygon
    sides=3, very thick, text width=1.5em },
  arn_e/.style = {treenode, blackblue, draw=blackblue,
    fill=bblue!10, circle,
    very thick, text width=1.5em },
  arn_w/.style = {treenode, black, draw=black,
    fill=white, circle,
    densely dashed, thick, text width=1.5em }
}

\newcommand{\tikzcircle}[2][red,fill=red]{\tikz[baseline=-0.5ex]\draw[#1,radius=#2] (0,0) circle ;}%

\usepackage{microtype}

\DeclareMathAlphabet\mathbfcal{OMS}{cmsy}{b}{n}

\newcommand{\mynewtheorem}[2]{
  \newaliascnt{#1}{dummy}
  \newtheorem{#1}[#1]{#2}
  \aliascntresetthe{#1}
  \expandafter\def\csname #1autorefname\endcsname{#2}
}

\theoremstyle{definition}
\mynewtheorem{theorem}{Theorem}
  \mynewtheorem{proposition}{Proposition}
  \mynewtheorem{corollary}{Corollary}
  \mynewtheorem{lemma}{Lemma}

\theoremstyle{definition}
  \mynewtheorem{remark}{Remark}
  \mynewtheorem{conjecture}{Conjecture}

\usepackage{etex}
\begin{document}

\title{How to generate random lambda terms?}
\author{Maciej Bendkowski${}^1$}
\address[1]{
      Theoretical Computer Science Department,\\
      Faculty of Mathematics and Computer Science,\\
      Jagiellonian University, {\L}ojasiewicza 6,\\
      30-348 Krak\'ow, Poland.
}
\email{maciej.bendkowski@\{tcs.uj.edu.pl,gmail.com\}}
\date{\today}
\maketitle

\begin{abstract}
We survey several methods of generating large random $\lambda$\nobreakdash-terms,
focusing on their closed and simply-typed variants. We discuss methods of exact- and
approximate-size generation, as well as methods of achieving size-uniform and
non-uniform outcome distributions.
\end{abstract}

\tableofcontents

\section{Introduction}\label{sec:introduction}
Generating examples of various combinatorial structures is an essential part of
devising counterexamples to working hypotheses as well as building up confidence
in their genuineness. It is a standard mathematical practice used to sieve out
and test promising ideas or conjectures before more rigorous and
labour-intensive treatment is called for. The more examples pass our tests, the
more evidence supports our hypothesis. To illustrate this point,
consider the famous $3n + 1$ Collatz conjecture.

\begin{conjecture}[Collatz, 1937]
  Let $f \colon \mathbb{N} \to \mathbb{N}$ be a function defined as
  \begin{equation}
    f(n) = \begin{cases}
      n / 2 & \text{if } n \equiv 0 \pmod{2}\\
      3n + 1 & \text{if } n \equiv 1 \pmod{2}
    \end{cases}
  \end{equation}
  Then, for all $n \geq 1$ there exists a sufficiently large $k$ such that
  $f^{(k)}(n) = 1$.
\end{conjecture}

The Collatz conjecture remains one of the most prominent open problems in
mathematics. It states that given \emph{any} natural number $n \geq 1$, the
infinite sequence
\begin{equation}
  f(n),~f(f(n)),~f(f(f(n))), \ldots
\end{equation}
must eventually reach the value $1$ or, equivalently, fall into the cycle $4
\to 2 \to 1 \to 4 \to \cdots$.

A myriad of computer-generated empirical evidence shows that it holds for all
positive integers up to $20 \times 2^{58} \approx 5.7646 \times
10^{18}$~\cite{lagarias2011}. Although no \emph{proof} of Collatz's conjecture
is known, the sheer bulk of empirical data provides high confidence in its
genuineness which, in turn, inspires new attempts at solving its mystery.

Nevertheless, even the highest levels of confidence cannot replace rigorous
mathematical reasoning. After all, we cannot test a hypothesis for all,
infinitely many different inputs. By a twist of perspective, however, all it
takes to \emph{disprove} a conjecture is a \emph{single} counterexample. That,
of course, can be looked for! A famous example of a conjecture which was refuted
by a single counterexample is Euler's sum of powers conjecture.

\begin{conjecture}[Euler, 1769]
  Let $n, k \geq 2$. If
    $a_1^k + a_2^k + \cdots + a^k_n = a^k_{n+1}$,
  then $n \geq k$.
\end{conjecture}

This long-standing conjecture was finally settled in the negative by Lander and
Parkin~\cite{lander1966}. Using a systematic, computer-assisted search procedure
they found the following neat counterexample to Euler's claim:
\begin{equation}
  {27}^5 + {84}^5 + {110}^5 + {133}^5 = {144}^5.
\end{equation} 
Without the use of computers, such a counterexample might not have been found.

\subsection{Why do we want random lambda terms?}
The \emph{generate-and-check} principle of testing working hypotheses is
mirrored in software development by \emph{software verification}. Since
\emph{proving} properties about large, industrial-strength programs is
tremendously difficult and time-consuming, in fact almost infeasible for most
types of applications, the modern practice is to \emph{test} them using a large
body of carefully crafted tests cases \emph{checking} the intended behaviour of
the program. Absolute, mathematical assurances are often given up in favour of
feasible, yet reasonable levels of \emph{confidence} in the actual program's
\emph{properties}.

Remarkably, test cases need not to be created manually by programmers but, can
(at least to some level) be generated \emph{automatically} by a computer. The
perhaps most notable examples of such a general testing process in software
development include model checking~\cite{cgjlv2000} and property-based software
testing such as QuickCheck~\cite{Claessen:2000:QLT}. In testing environments
like QuickCheck, the programmer is spared the burden of writing single test
cases. Instead, she is given the task of defining \emph{properties}
(\emph{invariants}) of tested software and providing a recipe for generating
random instances of input data. The testing environment generates hundreds of
random test cases and attempts to \emph{disprove} programmer-specified
properties. If a counterexample is found, the programmer is given
\emph{constructive evidence} that the indented property does not hold.
Otherwise, with hundreds or thousands of successful test cases, she can become
quite confident in the software's correctness.

One prominent application of random $\lambda$\nobreakdash-terms comes from
testing the correctness of certain program transformations in \textsc{ghc} ---
Haskell's most famous optimising compiler~\cite{marlow2012the}. Lambda terms are
a simple (if not \emph{the} simplest) example of a non-trivial combinatorial
structure with variables and corresponding name scopes. Consequently,
$\lambda$\nobreakdash-terms provide the scaffold of simple functional programs
and can be used to generate random inputs for compilers to consume. If, for some
specific input program $P$, the tested compiler generates two \emph{observably
different} executables (depending on the level of enabled optimisations) $P$
forms a constructive \emph{counterexample} to the working assumption that
performed optimisations do not change the semantics of the input program.
Remarkably, such a strikingly simple testing method lead to the discovery of
several important bugs in \textsc{ghc}'s strictness
optimiser~\cite{Palka:2011:TOC:1982595.1982615} (see
also~\cite{10.1145/3110259}). In particular, its misuse of the
\begin{Verbatim}[commandchars=\\\{\}]
  \PYG{n+nf}{seq} \PYG{o+ow}{::} \PYG{n}{a} \PYG{o+ow}{\PYGZhy{}\PYGZgt{}} \PYG{n}{b} \PYG{o+ow}{\PYGZhy{}\PYGZgt{}} \PYG{n}{b}
\end{Verbatim}
function of the following semantics:

\begin{Verbatim}[commandchars=\\\{\},codes={\catcode`\$=3\catcode`\^=7\catcode`\_=8}]
  \PYG{esc}{$\bot$} \PYG{p}{`}\PYG{n}{seq}\PYG{p}{`} \PYG{n}{b} \PYG{o+ow}{=} \PYG{esc}{$\bot$}
  \PYG{n}{a} \PYG{p}{`}\PYG{n}{seq}\PYG{p}{`} \PYG{n}{b} \PYG{o+ow}{=} \PYG{n}{b}
\end{Verbatim}

In the current paper we discuss several techniques of generating large, random
$\lambda$\nobreakdash-terms. We do not touch on the conceptually simpler
\emph{exhaustive generation} problem in which we are interested in the
construction of \emph{all} terms of a fixed size $n$, like in frameworks such as
SmallCheck~\cite{conf/haskell/RuncimanNL08}. Instead, we discuss the various
methods of \emph{sampling} (uniformly) random $\lambda$\nobreakdash-terms of
size $n$, including the two important classes of simply-typed and closed terms.

\subsection{What lambda terms do we want to sample?}
Before we begin to discuss \emph{how} to generate random
$\lambda$\nobreakdash-terms, let us pause for a moment and discuss \emph{what}
terms we want to sample. In order to sample (uniformly) random terms of some
size $n$, we have to ensure that there exists only a \emph{finite} number of
terms of that size. It means that we have to decide what to do with
$\alpha$\nobreakdash-equivalent terms and how to \emph{represent}
$\lambda$\nobreakdash-terms. Let us consider these issues in order.

Let $N$ and $M$ be two distinct $\alpha$\nobreakdash-equivalent
$\lambda$\nobreakdash-terms, differing only in names of their bound variables,
for instance $\lambda x y z. x z (y z)$ and $\lambda a b c. a c (b c)$. Although
syntactically different, both represent the same, famous $\mathbf{S}$
combinator. It can be therefore argued that under reasonable\footnote{Size
models in which the terms size does not depend on the particular names of bound
variables.} size models, both should be assigned the same \emph{size}. Since we
have an infinite supply of variable names, we are consequently inclined to
consider $\alpha$\nobreakdash-equivalent terms \emph{indistinguishable}. In
other words, focus on sampling $\alpha$\nobreakdash-equivalence classes of
$\lambda$\nobreakdash-terms, instead of their inhabitants.

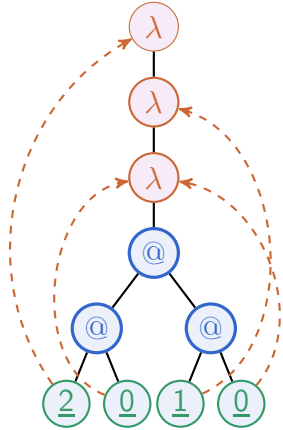
\begin{figure}[h]
  \centering
\begin{tikzpicture}[>=stealth',level/.style={thick}]
\draw
node[arn_r](A){\( \lambda \)}
child[level distance=1cm]{
node[arn_r](B){\( \lambda \)}
child[level distance=1cm]{
node[arn_r](C){\( \lambda \)}
child[level distance=1cm]{
    node[arn_n]{ \( @ \) }
    child[level distance=1cm]{
        node[arn_n]{ \( @ \) }
        child[sibling distance=.8cm]{
            node[arn_g](D){ $\idx{2}$ }
        }
        child[sibling distance=.8cm]{
            node[arn_g](E){ $\idx{0}$ }
        }
    }
    child[level distance=1cm]{
        node[arn_n]{ \( @ \) }
        child[sibling distance=.8cm]{
            node[arn_g](F){ $\idx{1}$ }
        }
        child[sibling distance=.8cm]{
            node[arn_g](G){ $\idx{0}$ }
        }
    }
}
}
};
\path[->] (D) edge [bred,dashed,thick,bend left=48] node {} (A);
\path[->] (E) edge [bred,dashed,thick,bend left=75] node {} (C);
\path[->] (F) edge [bred,dashed,thick,bend right=70] node {} (B);
\path[->] (G) edge [bred,dashed,thick,bend right=60] node {} (C);
\end{tikzpicture}
\caption{Tree-like representation of the term $\lambda \lambda \lambda . \idx{2}
\idx{0} (\idx{1} \idx{0})$. Note that the traditionally implicit term
application is presented in an explicit form.}\label{fig:example:term}
\end{figure}

Given these concerns, we have to decide on a representation in which we do not
differentiate bound variable names. We can therefore \emph{pick} one
\emph{canonical} representative for each $\alpha$\nobreakdash-equivalence class
of terms, in effect sampling terms \emph{up to}
$\alpha$\nobreakdash-equivalence,
see~\cite{Wang05generatingrandom,Wang05generatingrandom2,BODINI2013227}, or use
a representation in which there exists just one such representative, for
instance the De~Bruijn notation~\cite{deBruijn1972},
see~\cite{10.5555/2438099.2438158,BODINI201845,bendkowski_grygiel_tarau_2018}.
Recall that in the De~Bruijn notation there are no named variables. Instead, we
use \emph{dummy} indices $\idx{0},\idx{1},\idx{2},\ldots$ to denote variable
occurrences. The index $\idx{n}$ represents a variable which is bound by its
$(n+1)$st proceeding abstraction. Should there be fewer abstractions than $n+1$,
then $\idx{n}$ represents a free variable occurrence. And so, in the De~Bruijn
notation both $\lambda x y z. x z (y z)$ and $\lambda a b c. a c (b c)$ admit
the same representation $\lambda \lambda \lambda . \idx{2} \idx{0} (\idx{1}
\idx{0})$, see~\Cref{fig:example:term}.

Neither representation is better than the other one. What makes them quite
different, however, are the \emph{details} of related size models, in particular
what \emph{weights} we assign to abstractions, applications, and variables. For
instance, how should measure the size of a variable? Should it be proportional
to the distance between the variable occurrence and its binding abstraction, or
should it be independent of that distance? What toll should be impose on the
binder distance? Should it be a linear function or perhaps a more sophisticated
one, say logarithmic? If we use De~Bruijn indices, how should we represent them?
Should we use the more traditional unary base where $\idx{n}$ is just an
$n$\nobreakdash-fold application of the successor function to zero, or should we
represent indices in binary format?

Naturally, these representation details prompt \emph{different} size models and,
consequently, alter the landscape of random $\lambda$\nobreakdash-terms we are
going to sample. What details should we therefore agree on? Depending on the
specifics of the representation and size model, large random
$\lambda$\nobreakdash-terms change their \emph{statistical characteristics}.
Knowing these traits, we can choose a specific size model (and representation)
which best suits our needs. For instance, large random
$\lambda$\nobreakdash-terms in representations in which variables have constant
\emph{weight} tend to be \emph{strongly-normalising}~\cite{lmcs:848}. On the
other hand, in De~Bruijn models in which variables have size proportional to the
distance between them and their binding abstractions, large random terms are
almost never \emph{strongly-normalising}~\cite{10.1093/logcom/exx018}.

Finding properties of large, random $\lambda$\nobreakdash-terms is a fascinating
and surprisingly challenging endeavour. It often involves sophisticated counting
arguments using generating functions and advanced methods of analytic
combinatorics~\cite{flajolet09}. For some size models, however, even these
powerful techniques are not so easily applicable. In the current survey, we do
not dive into theoretical details of these methods. Instead, we focus more on
their \emph{application} to the random generation of
$\lambda$\nobreakdash-terms. We invite the curious reader to follow the
referenced literature for more details.

\begin{remark}
 We are interested in generating large, uniformly random
 $\lambda$\nobreakdash-terms. With large, \emph{unbiased} terms it is
 possible to test not only the \emph{correctness} of a compiler or abstract
 machine implementing the target language, but also test its \emph{performance} and
 \emph{scalability}, otherwise impossible to test through small input programs.

 Our decision has a few notable consequences. We are \emph{not} interested in
type-oriented generation methods similar to the Ben-Yelles
algorithm~\cite{benyelles} or ad-hoc techniques which do not give the user
direct and rigorous control over the outcome term distribution. While we are
primarily interested in \emph{uniform} distributions, we also address the issue
of \emph{non-uniform generation} in~\Cref{sec:multiparametric:samplers}.
\end{remark}

\section{Combinatorial generation methods}
Generating random combinatorial structures is a well-established topic in
theoretical computer science. Let us examine how standard generation methods can be tailored
to construct uniformly random $\lambda$\nobreakdash-terms.

\subsection{Exact-size sampling and the recursive
method}\label{subsec:recursive:method} Nijenhuis and Wilf’s recursive method
\cite{NijenhuisWilf1978}, later systematised by Flajolet, Zimmerman and Van
Cutsem~\cite{FLAJOLET19941}, is a simple, general-purpose recursive sampling
template for a wide range of combinatorial objects, including context-free
languages and various tree-like structures.

Given a formal specification of so-called \emph{decomposable} structures, the
recursive method provides an inductive scheme of constructing recursive samplers
meant to generate objects of some fixed size $n$.  In what follows we
adopt the recursive method and design a simple sampler for closed
$\lambda$\nobreakdash-terms in the De~Bruijn notation:
\begin{equation}\label{eq:combinatorial:generation:lambda:terms:def}
    N, M := \idx{n}~|~(N M)~|~\lambda N.
\end{equation}

Following the recursive method, terms will be built in a top-down fashion,
according to their inductive
specification~\eqref{eq:combinatorial:generation:lambda:terms:def}. In fact, we
are going to devise a slightly more general recursive sampler
$\textsc{Gen}_{\Gamma}(n)$ building terms of size $n$, with free indices in the
set $\Gamma$. The desired closed $\lambda$\nobreakdash-term sampler
$\textsc{Gen}_{\emptyset}(n)$ will be recovered afterwards.

During its execution, $\textsc{Gen}_{\Gamma}(n)$ will make a series of random
choices. Each of these random decisions will we consulted with an external
oracle $\CS{O}$ providing suitable \emph{branching probabilities} ensuring a
uniform outcome distribution. We explain later how to construct such an oracle
given the formal $\lambda$\nobreakdash-term
specification~\eqref{eq:combinatorial:generation:lambda:terms:def}.

When invoked, $\textsc{Gen}_{\Gamma}(n)$ has to decide whether to output one of
the available indices $\idx{n}\in \Gamma$, an application of two terms $(N M)$,
or an abstraction $\lambda N$. The sampler queries the oracle $\CS{O}$ for
respective branching probabilities, and chooses a constructor according to the
obtained distribution. Building an index $\idx{n} \in \Gamma$ is trivial.
Assuming that we have enough \emph{fuel} $n$ to pay for the index $\idx{n}$, we
just need to look it up in $\Gamma$. Applications and abstractions are a bit
more involved.

Note that if $\lambda N$ is a term with free indices in $\Gamma$, then the free
indices of its subterm $N$ are elements of $\Gamma' = \Gamma_\uparrow \cup
\{\idx{0}\}$ where $\Gamma_\uparrow = \{ \idx{n+1}\colon \idx{n} \in \Gamma\}$
is the \emph{lifted} variant of $\Gamma$.  Therefore, if
$\textsc{Gen}_{\Gamma}(n)$ decides to construct an abstraction $\lambda N$, it
can build $N$ by invoking $\textsc{Gen}_{\Gamma'}(n - a)$, accounting $a$ size
units for the head abstraction (depending, of course, on the assumed size model).

Let us consider what happens when $\textsc{Gen}_{\Gamma}(n)$ decides
to generate an application $(N M)$.  Denote the size contribution of a single
application as $b$. Now, in order to build $(N M)$ we first query the oracle
$\CS{O}$ for a random $i=1,\ldots,n-1$. Next, we construct two terms $N,M$ of sizes $i$
and $n - i - b$, respectively. Note that variable contexts of both $N$ and $M$
do not change. We can therefore readily invoke $\textsc{Gen}_{\Gamma}(i)$ and
$\textsc{Gen}_{\Gamma}(n -i - b)$, and apply $N$ to $M$.

\medskip
{\bf Oracle construction.}
In order to make its random decisions $\textsc{Gen}_{\Gamma}(n)$ repeatedly
queries an external oracle $\CS{O}$. Suppose that we invoke
$\textsc{Gen}_{\Gamma}(n)$. How should $\CS{O}$ compute the necessary branching
probabilities? The idea is quite simple.

Let $L_{\Gamma}(n)$ denote the set of $\lambda$\nobreakdash-terms of size $n$
with free indices in the set $\Gamma$. In order to assign each term in
$L_\Gamma(n)$ a uniform probability
\begin{equation}
    p = \frac{1}{|L_\Gamma(n)|}
\end{equation}
we first determine $|L_\Gamma(n)|$. Then, the branching probabilities $p_1, p_2$
and $p_3$ corresponding to choosing and index $\idx{n}$, an abstraction $\lambda
N$, and an application $(N M)$, respectively, are given by
\begin{eqnarray}
  \begin{split}
  p_1 &= p \cdot |\{\idx{n} \colon \idx{n} \in L_{\Gamma}(n)\}|\\
  p_2 &= p \cdot |\{\lambda N \colon \lambda N \in
  L_{\Gamma}(n)\}|\\
  p_3 &= p \cdot |\{(N M) \colon (N M) \in L_{\Gamma}(n)\}|.
  \end{split}
\end{eqnarray}

Likewise, the probability of choosing an application $(N M)$ in which $N$ is of
size $i$ and $M$ is of size $n-i-b$ is equal to the proportion of such
applications among all applications of size $n$, all in the common context
$\Gamma$. It means that in order to compute appropriate branching probabilities
our oracle needs to calculate cardinalities of involved sets.

Computing the cardinality of $L_{\Gamma}(n)$, as well as specific subsets can be
achieved using dynamic programming, much like $\textsc{Gen}_{\Gamma}(n)$ itself.
The oracle $\CS{O}$ can be precomputed once, memorised, and reused in between
subsequent sampler invocations.

\medskip
{\bf Complexity and limitations.}
The recursive method, albeit simple, admits considerable practical limitations.
Both the oracle construction and the following sampling procedure require
$\Theta(n^2)$ arithmetic operations on integers exponential in the target size
$n$, turning the effective bit complexity of the recursive method to
$O(n^{3+o(1)})$. It is therefore possible to sample
$\lambda$\nobreakdash-terms of moderate sizes in the thousands.

Let us remark that Denise and Zimmerman combined the recursive method with
certified floating-point arithmetic, reducing the expected time and space
complexity of sampling down to $O(n^{1+o(1)})$ and $O(n^{2+o(1)})$ preprocessing
time~\cite{DenZimm99}. Notably, for context-free languages the preprocessing
time can be further reduced to $O(n^{1 + o(1)})$.

\begin{remark}\label{sec:random:generation:remark:lescanne:grygiel}
    Grygiel and Lescanne~\cite{grygiel_lescanne_2013} proposed a somewhat
    similar sampler for closed $\lambda$\nobreakdash-terms based on the
    following \emph{ranking-unranking} principle.

    Given a set $S$ of $n$ elements we wish to sample from, we construct a pair
    of mutually inverse bijections $f \colon S \to \{0,1,\ldots,n-1\}$ and
    $f^{-1} \colon \{0,1,\ldots,n-1\} \to S$. The former function defines the
    \emph{rank} $f(s)$ of an element $s \in S$, whereas the latter
    \emph{unranking} function maps a given rank $f(s)$ to the unique element $s
    \in S$ attaining its value. Both $f$ and its inverse $f^{-1}$ are
    constructed in a recursive manner. With these two functions, sampling
    elements of $S$ boils down to sampling a uniformly random number $i =
    0,\ldots,n-1$. Once we obtain a random rank $i$, we can use the unranking
    function $f^{-1}$ to construct the corresponding element $f^{-1}(i)$.

    Based on the ranking-unranking sampler for closed
    $\lambda$\nobreakdash-terms (in the size model where variables do not
    contribute to the term size) Grygiel and Lescanne derived a rejection
    sampler for simply-typed $\lambda$\nobreakdash-terms.
    Such a sampler generates closed (untyped) terms until a typeable one is
    sampled. Intermediate terms are simply discarded.

    The efficiency of such a rejection scheme depends on two factors --- the
    efficiency of the underlying sampler for untyped terms and, even more importantly, the
    \emph{proportion} of typeable terms among closed (untyped) terms. Unfortunately,
    already for $n = 50$ the ratio between typeable terms and untyped ones is
    less than $10^{-5}$, see~\cite[Section 9.3]{grygiel_lescanne_2013}. It means
    that the average number of trials required to find a typeable term of size
    $50$ is already of order $10^5$. Such a sampler is not likely to
    output terms larger than $n = 50$.

    Let us remark that Lescanne devised a more direct sampler for \emph{linear}
and \emph{affine} $\lambda$\nobreakdash-terms based on the same
ranking-unranking principle~\cite{Lescanne:2018:QAL:3185755.3173547}. Tracking
the evolution of respective (simpler) contexts, Lescanne proposed samplers for a
few existing size models.
\end{remark}

\begin{remark}
    Wang proposed a recursive sampler for closed $\lambda$-terms in a size model
    where all constructors (\emph{i.e.}~variables, applications, and abstractions)
    contribute one to the overall term size~\cite{Wang05generatingrandom}.  She
    then adopted her samplers to simply-typed $\beta$\nobreakdash-normal
    forms~\cite{Wang05generatingrandom2} using a finite truncation of Takahashi,
    Akama, and Hirokawa's context-free grammars $G(\Gamma,\alpha)$ generating
    $\beta$\nobreakdash-normal forms of type $\alpha$ in the context
    $\Gamma$~\cite{TAKAHASHI1996144}.

    Since her sampler is based on a dedicated grammar for
    $\beta$\nobreakdash-normal forms, it cleverly avoids the slowdown imposed by
    rejection sampling from the larger set of untyped terms.  Unfortunately,
    such a method comes at a significant price.  Grammar-based samplers
    do not scale onto larger types. Involved grammars $G(\Gamma,\alpha)$
    quickly become intractable and too expensive to be effectively constructed.

    Let us remark that using the idea of truncated type grammars Asada, Kobayashi,
    Sin’ya and Tsukada showed recently that, asymptotically almost all typeable
    $\lambda$\nobreakdash-terms of order $k$ have a $(k-1)$-fold exponential
    reduction sequence~\cite{lmcs:5203}.
\end{remark}

\subsection{Generating functions and the analytic toolbox}
Counting and generating random objects are intimately connected activities. In order to construct
random $\lambda$\nobreakdash-terms we usually need to compute a series of branching
probabilities. These, in turn, depend on the specific cardinalities of many
different term families for various size values.

In order to conveniently operate on (infinite) \emph{counting sequences}
enumerating respective cardinalities, we resort to \emph{generating functions},
\emph{i.e.}~formal power series $\sum_{n \geq 0} f_n z^n$ whose coefficients
$\seq{f_n}_n$ encode the counting sequence we are interested in. To illustrate the convenience of
generating functions, consider the famous example of \emph{Catalan numbers}
$\seq{c_n}_n$ enumerating, \emph{inter alia}, the number of expressions
containing $n$ pairs of matched parentheses, or plane binary trees with $n$
internal nodes. These numbers satisfy the following
\emph{recursive} relation:
\begin{equation}\label{eq:catalan:sequence}
  c_{n+1} = \sum_{i = 0}^{n} c_i \cdot c_{n-i} \qquad \text{where} \qquad c_0 = 1.
\end{equation}

With the help of generating functions, its possible to represent the entire
infinite counting sequence of Catalan numbers as a \emph{finite} expression.
Lifting~\eqref{eq:catalan:sequence} onto the level of generating functions
involving $C(z) = \sum_{n \geq 0} c_n z^n$ we
find that
\begin{equation}\label{eq:catalan:gf}
  C(z) = 1 + z {C(z)}^2 = \frac{1-\sqrt{1-4z}}{2z}.
\end{equation}
In this form, we can recover the $n$th Catalan number by finding the
$n$th coefficient $[z^n]C(z) = C^{(n)}(0) / n!$ of $C(z)$ using the Taylor series
expansion of $C(z)$ around $z = 0$.

Typically, when solving recursive equations like~\eqref{eq:catalan:sequence} we
do not concern ourself with the convergence of associated generating functions.
If, however, the generating functions happen to be convergent, they correspond to
complex \emph{analytic functions}. We can then draw from the rich fountain of
analytic combinatorics~\cite{flajolet09} --- a theory devoted to giving
precise quantitative predictions of large combinatorial structures, exploiting
for that purpose the analytic properties of related generating functions.

\medskip
{\bf Holonomic functions and P-recursive sequences.}
We say that a formal power series ${f(z) = \sum_{n \geq 0} f_n z^n}$ is
\emph{holonomic
(D-finite)} if it satisfies a linear differential equation
\begin{equation}\label{eq:random:generation:holonomic}
    P_r(z) \frac{d^r}{d z^r} f(z) + P_{r-1}(z) \frac{d^{r-1}}{d z^{r-1}} f(z) +
    \cdots + P_1(z) \frac{d}{d z} f(z) + P_0(z) f(z) = 0
\end{equation}
for some polynomials $P_i(X) \in \mathbb{C}[X]$, \emph{cf.}~\cite[Appendix
B.4]{flajolet09}. By extension, a function $f \colon \mathbb{C} \to \mathbb{C}$
analytic at $0$ is said to be \emph{holonomic} if its power series representation is
holonomic. Likewise, the coefficient sequence $\seq{f_n}_n$ of a holonomic power
series $f(z)$ is called \emph{holonomic} (\emph{P-recursive}).

Holonomic functions are a rich subclass of analytic functions emerging naturally
in the context of various enumeration problems. Notably, D-finite functions
enjoy a series of pleasing properties. For instance, holonomic generating
functions subsume rational and algebraic functions --- two classes of generating
functions linked to counting
problems involving rational and context-free languages. Given an algebraic
function $f(z)$ represented as a branch of a polynomial equation $P(z, f(z)) =
0$ it is possible to compute its corresponding D-finite
representation~\cite{comtet}. In fact, modern computer algebra systems provide
dedicated packages specifically for that purpose,
\emph{cf.}~\cite{SalvyZimmermann1994}.

D-finite functions have found numerous applications in combinatorial enumeration
and symbolic computations~\cite{aequalsb}. Remarkably, using respective linear
differential equations~\eqref{eq:random:generation:holonomic} it is possible to
find linear recurrences defining their coefficient sequences $\seq{f_n}_n$ and,
consequently, compute the $n$th coefficient $f_n$ using just $O(n)$ arithmetic
operations.

Consider our running example of Catalan numbers.
The related holonomic equation for $C(z)$ takes the form
\begin{equation}\label{eq:random:generation:holonomic:2}
    1+ \left( 2 z - 1\right) C(z) + \left( 4{z}^{2}-z
    \right) {\frac {d}{{d}z}} C(z) = 0.
\end{equation}

In the special case of generating functions, the differential operator
$\frac{d}{{d}z}$ admits a natural \emph{combinatorial} interpretation. Note that
the defining relation
\begin{equation}\label{eq:pointing:def}
  z \frac{d}{{d}z} C(z) = z \frac{d}{{d}z}\sum_{n \geq 0} c_n z^n = \sum_{n \geq
1} n \cdot c_n z^{n}
\end{equation}
implies that the coefficient $[z^n] z \frac{d}{d z}C(z)$ is equal to $n \cdot
[z^{n}]C(z)$. We can therefore interpret the generating function $z
\frac{d}{{d}z}C(z)$ as encoding the counting sequence of plane binary trees with
\emph{pointed} internal nodes. Each tree of size $n$ has $n$ internal nodes and
therefore gives rise to $n$ different variants of itself, each with a different
node being pointed.

Using~\eqref{eq:pointing:def} we can translate the holonomic equation defining
$C(z)$ into a linear recurrence involving its coefficients:
\begin{equation}\label{eq:random:generation:holonomic:3}
    \left(n + 2\right) c_{n+1} = 2 \left( 2\,n + 1 \right) c_n \quad
    \text{with} \quad c_0 = 1.
\end{equation}

It is clear that with the help of~\eqref{eq:random:generation:holonomic:3} we
can compute the $n$th Catalan number $c_n$ using only a linear number of
arithmetic operations. On the other hand, using the naïve method derived from
the defining identity $C(z) = 1 + z {C(z)}^2$, we have to compute a sum of
products of all $c_0,c_1,\ldots,c_n$ in order to compute the next Catalan number
$c_{n+1}$. Such a process requires $\Theta(n^2)$ arithmetic operations.

Following similar steps, each holonomic (in particular algebraic) function
admits a linear recurrence governing its coefficients. Note that by exploiting
this fact we can, for instance, speed up the oracle computations involved in the
recursive method.

\subsection{Combinatorial bijections and Rémy's algorithm}
The simple recurrence~\eqref{eq:random:generation:holonomic:3} defining $c_n$
provides the basis for an elegant sampling algorithm for plane binary trees due
to Rémy~\cite{Remy85}. The main idea is to interpret it \emph{combinatorially},
using it as an indication on how to map the set of binary trees with $n$
internal nodes to the larger set of trees with $n+1$ internal nodes.

We start with noticing that each binary tree with $n$ internal
nodes has $n+1$ leaves and so the total of $2n+1$ nodes. It means that we can
interpret~\eqref{eq:random:generation:holonomic:3} as a bijection between two
sets of trees --- on the left-hand side, the set of binary trees with $n+1$
internal nodes, each with a single \emph{pointed} leave; and on the right-hand
side, the set of binary trees with $n$ internal nodes, each with a single pointed node
coloured with one of two distinct colours, say red \tikzcircle[fill=red]{3pt}
or blue \tikzcircle[fill=bblue]{3pt}.
Let $T$ be a binary tree with $n$ internal nodes.  Consider
the following \emph{grafting} operation:
\begin{itemize}
    \item Select a random node (be it internal or external) in $T$ and call it $x$.
    \item Flip a coin. Depending on the outcome, colour $x$ either red~\tikzcircle[fill=red]{3pt}
or blue~\tikzcircle[fill=bblue]{3pt}.
\item If $x$ was coloured red \tikzcircle[fill=red]{3pt} replace it with \raisebox{-1.0 ex}{
\begin{tikzpicture}[level/.style={sibling distance = 0.4cm/#1, level
distance = .3cm}]
\draw
++(0, -1.5)
    node[arn_nn]{\( \)}
    child{
        node[arn_rp]{\( { \scriptscriptstyle x } \)}
    }
    child{
        node[arn_nn]{\( \)}
    }
;
\end{tikzpicture}
} in $T$.
\item Symmetrically, if $x$ was coloured blue \tikzcircle[fill=bblue]{3pt} replace it with \raisebox{-1.0 ex}{
\begin{tikzpicture}[level/.style={sibling distance = 0.4cm/#1, level
distance = .3cm}]
\draw
++(0, -1.5)
    node[arn_nn]{\( \)}
    child{
        node[arn_nn]{\( \)}
    }
    child{
        node[arn_rp]{\( { \scriptscriptstyle x } \)}
    }
;
\end{tikzpicture}
} in $T$.
\item Point to the newly created leave (sibling of $x$).
\end{itemize}

Suppose that we started the grafting operation with a (uniformly) random tree
with $n$ internal nodes. Note that once we select a random node $x$ and assign
it a colour, we end up with a random tree $T$ enumerated by the right-hand side
of~\eqref{eq:random:generation:holonomic:3}. At this point we readily notice
that after the above grafting operation we obtain a random tree enumerated by
the left-hand side of~\eqref{eq:random:generation:holonomic:3}. We can therefore
forget the pointer created in the final grafting step and obtain a random binary
tree with $n+1$ internal nodes.

Clearly, this process can be iterated as many times as required. If we want to
sample a random binary tree of size $2n + 1$, we can start with a single node
\begin{tikzpicture}[level/.style={sibling distance = 0.4cm/#1, level distance =
.3cm}] \draw ++(0, -1.5) node[arn_nn]{\( \)}; \end{tikzpicture} and repeat the
grafting operation $n - 1$ times. In his celebrated \emph{Art of Computer
Programming}~\cite{Knuth:2006:ACP:1121689} Knuth provides an elegant
implementation of Rémy's algorithm which he calls algorithm \textbf{R}. With its
help, we obtain an efficient, linear time, exact-size sampling algorithm for
plane binary trees.

\begin{remark}
     Finding constructive interpretations for combinatorial identities is, in
     general, a non-trivial and creative task.  Rémy's algorithm is somewhat
     \emph{ad-hoc} and therefore cannot be easily generalised onto other tree
     structures. Interestingly enough, let us remark that Bacher, Bodini, and
     Jacquot were able to use ideas of Boltzmann sampling and exploit the
     holonomic specification for so-called Motzkin trees (\emph{i.e.}~unary-binary
     trees), developing a linear time sampler~\cite{BBJ}.
\end{remark}

\begin{remark}
Quite often appropriate bijections are derived only after combinatorial
identities involving respective generating function are established. For
instance, let us mention the non-trivial bijections between combinatorial maps
and certain classes of enriched trees leading to efficient samplers for linear
and affine $\lambda$\nobreakdash-terms in the \emph{canonical} representation of
$\lambda$\nobreakdash-terms constructed up to
$\alpha$\nobreakdash-equivalence~\cite{BODINI2013227}, or the bijection between
neutral $\lambda$\nobreakdash-terms and Motzkin trees in the so-called
\emph{natural} size model under the De~Bruijn
notation~\cite{10.1093/logcom/exx018}.
\end{remark}

\medskip
{\bf Random combinators.}
Rémy's algorithm is an important sampling algorithm allowing us to generate
combinators of any finite bases of primitive combinators. To illustrate this
point, consider the example of $\mathbf{S K}$\nobreakdash-combinators specified
as
\begin{equation}\label{eq:random:generation:combinators}
    C := \mathbf{S}~|~\mathbf{K}~|~(C C).
\end{equation}

Assume that the size of a combinator is equal to the number of its internal
applications\footnote{If $\mathbf{S}$ and $\mathbf{K}$ contribute to the size,
the corresponding counting sequence becomes periodic and introduces some
technical (though manageable) difficulties, similarly to the case of
leaves in plane binary trees.}. Suppose that we want to sample a combinator of
size $n$. Note that we can \emph{decompose} each combinator of size $n$ into a
plane binary scaffold with $n$ internal nodes, and a sequence of $n+1$ primitive
combinators $\mathbf{S}$ and $\mathbf{K}$.  The scaffold determines the
application structure of the combinator whereas the sequence governs the
placement of $\mathbf{S}$ and $\mathbf{K}$s in the term.  Moreover, such a
decomposition is clearly \emph{invertible} --- simply traverse the scaffold in-order
and assign elements of the sequence to successive leaves.

The unambiguous decomposition of combinators allows us to use Rémy's algorithm
to sample a random scaffold with $n$ internal nodes and, independently, generate
a sequence of $n+1$ random primitive combinators. Afterwards, we compose the two
into a uniformly random combinator.

\begin{remark}
The above sampling scheme is quite efficient. It allows us to easily sample random
    combinators even of sizes in the hundreds of millions. Unfortunately, even with
    the Motzkin tree sampler of Bacher, Bodini, and Jacquot~\cite{BBJ} it is not so easy
    to provide an efficient, analogous sampler for closed $\lambda$\nobreakdash-terms.

    Unlike binary trees, Motzkin trees of some fixed size $n$ might have a
    varying number of leaves. Moreover, in order to interpret a leave as a
    variable $x$, it is important to know the number of its potential binders,
    \emph{i.e.}~unary nodes between $x$ and the term root. We refer the curious reader
    interested in empirical evaluation of such scaffold decompositions for
    $\lambda$\nobreakdash-terms
    to~\cite{10.1007/978-3-319-73305-0_8,10.1007/978-3-319-94460-9_15}.
\end{remark}

\subsection{Boltzmann models}
For many years, the exact-size sampling paradigm was the \emph{de facto}
standard in combinatorial generation. Its long-lasting dominance was brought to
an end with the seminal paper of Duchon, Flajolet, Louchard, and
Schaeffer~\cite{Duchon:2004:BSR:1024662.1024669} who introduced the framework of
\emph{Boltzmann models}. Instead of generating random structures of fixed size
$n$, Boltzmann models provide a more relaxed, \emph{approximate-size} sampling
environment which uses the analytic properties of associated generating functions.

Fix a class $\CS{A}$ of combinatorial structures we want to sample.  Under the
\emph{Boltzmann model} we assign to each object $\alpha \in \CS{A}$ a
probability measure
\begin{equation}\label{eq:random:generation:boltzmann:prob}
    \mathbb{P}_x(\alpha) = \frac{x^{|\alpha|}}{A(x)}
\end{equation}
where $x$ is some aptly chosen, positive real-valued \emph{control parameter},
and $|\alpha|$ denotes the size of object $\alpha$. If we sample in accordance
with the above Boltzmann distribution, outcome structures of equal size
invariably occur with \emph{the same} probability. Boltzmann models retain
therefore the usual requirement of \emph{uniformity} however now, the size of
the outcome object is no longer fixed, but instead varies. Indeed, let $N$ be a
random variable denoting the outcome size of a \emph{Boltzmann sampler}
generating objects according to~\eqref{eq:random:generation:boltzmann:prob}. Summing over all
$\alpha \in \CS{A}$ of size $n$ we find that
\begin{equation}
    \mathbb{P}_x(N = n) = \frac{a_n x^n}{A(x)}.
\end{equation}

We can control the expected size or generated objects by choosing suitable
values of the control parameter $x$. The expected average and standard
deviation of $N$ satisfy, respectively,
\begin{equation}
    \mathbb{E}_x(N) = x \frac{A'(x)}{A(x)} \quad \text{and} \quad
    \sigma_x(N) = \sqrt{x \mathbb{E}'(N)}.
\end{equation}

Suppose, for instance, that we want to design a Boltzmann model for plane binary
trees. Recall that the associated generating function $C(z)$ satisfies the
relation $C(z) = 1 + z {C(z)}^2$. If we want to centre the mean output tree size
around some fixed value $n$, we need to find a suitable control parameter $x$
such that $\mathbb{E}_x(N) = n$. A direct computations reveals that we should
use
\begin{equation}\label{eq:random:generation:boltzmann:2}
x = \frac{n (n+1)}{(2 n+1)^2}.
\end{equation}

So, if we want to obtain a random tree in some admissible size window
${[(1-\varepsilon)n, (1+\varepsilon)n]}$ with \emph{tolerance} $\varepsilon$, we
calibrate $x$ according to~\eqref{eq:random:generation:boltzmann:2} and reject
trees falling outside of the prescribed size window.

\medskip
{\bf Boltzmann samplers.}
For many interesting instances of decomposable combinatorial classes, such as
rational languages or algebraic data types, it is possible to
\emph{automatically} design appropriate Boltzmann models and corresponding
samplers generating random objects in accordance with the underlying Boltzmann
distribution,
see~\cite{Duchon:2004:BSR:1024662.1024669,Canou:2009:FSR:1596627.1596637,doi:10.1137/1.9781611975062.9,PIVOTEAU20121711}.

Much like exact-size recursive samplers, the Boltzmann sampler layout follows
the inductive structure of the target specification.  Suppose that we want to
sample an object from a (disjoint) union class $\CS{A} = \CS{B} + \CS{C}$. According to the
Boltzmann distribution~\eqref{eq:random:generation:boltzmann:prob} the
probability of sampling an object from $\CS{B}$ is equal to $B(x) / A(x)$.
Likewise, the probability of sampling an object from $\CS{C}$ is equal to $C(x)
/ A(x)$. The corresponding Boltzmann sampler $\Gamma(\CS{A})$ for $\CS{A}$
performs therefore a single Bernoulli trial with parameter $B(x) / A(x)$.  If
successful, the sampler invokes $\Gamma(\CS{B})$. Otherwise, it calls the
remaining sampler $\Gamma(\CS{C})$. Clearly, such a sampler conforms with the
Boltzmann distribution.

Now, suppose that we want to sample an object from a product class $\CS{A} =
\CS{B} \times \CS{C}$. According to the Boltzmann distribution, the probability
measure of a \emph{product object} $\alpha = (\beta,\gamma) \in \CS{A}$ satisfies
\begin{equation}
  \mathbb{P}_x(\alpha) = \frac{x^{|\alpha|}}{A(x)} = \frac{x^{|\beta| + |\gamma|}}{A(x)} = \frac{x^{|\beta|}
    x^{|\gamma|}}{A(x)} = \mathbb{P}_x(\beta) \cdot \mathbb{P}_x(\gamma).
\end{equation}

Therefore, in order to sample an object from $\CS{A}$ we can independently
invoke samplers $\Gamma(\CS{B})$ and $\Gamma(\CS{C})$, collect their outcomes,
and construct a pair out of them. The result is a random object from $\CS{A}$.

\begin{remark}
The disjoint union $+$ and Cartesian product $\times$ suffice to design
    Boltzmann samplers for algebraic specifications, in particular unambiguous
    context-free grammars. Let us remark, however, that Boltzmann sampler are
    not just limited to these two operations.
    Over the years, efficient Boltzmann samplers for both labelled and unlabelled
    structures were developed. Let us just mention Boltzmann samplers for
    various Pólya structures~\cite{Flajolet:2007:BSU:2791135.2791140}, first-order
    differential specifications~\cite{BODINI20122563}, labelled planar
    graphs~\cite{fusy2005quadratic}, or plane partitions~\cite{bodini2010random}.

    In most cases, especially when optimised using \emph{anticipated
rejection}~\cite{BodGenRo2015}, Boltzmann samplers are remarkably efficient,
frequently generating random structures of sizes in the millions. For typical
algebraic specifications, Boltzmann samplers calibrated to the target size
window $[(1-\varepsilon)n, (1+\varepsilon)n]$ admit, on average, a linear $O(n)$
time complexity under the real arithmetic computation
model~\cite{Duchon:2004:BSR:1024662.1024669,BodGenRo2015}. For exact-size
sampling where the tolerance $\varepsilon = 0$, the complexity becomes $O(n^2)$.
\end{remark}

\begin{remark}
Boltzmann samplers for $\lambda$\nobreakdash-terms were first developed by
    Lescanne~\cite{DBLP:journals/corr/Lescanne14,grygiel_lescanne_2015} for
    Tromp's binary $\lambda$\nobreakdash-calculus~\cite{tromp:DSP:2006:628}.
    Tromp's special encoding of (plain) lambda terms as binary strings, provided
    an unambiguous context-free specification fitting the framework of
    Boltzmann samplers. With their help, it became possible to generate
    $\lambda$\nobreakdash-terms a couple of orders of magnitude larger than using
    recursive method.

    Boltzmann models, especially if used with rejection methods, provide a rich
source of random $\lambda$\nobreakdash-terms. Although quite general, Boltzmann
samplers cannot be used for all term representations and size notions. Recall
that in order to derive a corresponding Boltzmann model, we have to, \emph{inter
alia}, evaluate the associated generating functions at the calibration parameter
$x$. For some representations, such as for instance Wang's closed
$\lambda$\nobreakdash-term model~\cite{Wang05generatingrandom}, the associated
generating function is not analytic around the complex plane origin\footnote{In
such a case the corresponding counting sequence grows super-exponentially fast,
\emph{i.e.}~asymptotically faster than any exponential function $C^n$.}. Moreover,
without a finite specification we cannot derive a finite set of Boltzmann
samplers. In such cases we have to resort to alternative sampling methods
including, \emph{e.g.}~rejection-based samplers.
\end{remark}

\subsection{Rejection samplers}
Most interesting subclasses of $\lambda$\nobreakdash-terms, such as closed or
(simply) typeable terms, have no known finite, \emph{admissible}\footnote{Purely
  in
  terms of basic constructions like the disjoint sum or Cartesian product,
  cf.~\cite[Part A. Symbolic Methods]{flajolet09}.} specification which
would allow us to derive effective Boltzmann samplers. In order to sample from
such classes, we \emph{approximate} them using finite, admissible specifications,
using rejection methods if needed.

\medskip
{\bf Closed $\lambda$-terms.}
Most $\lambda$\nobreakdash-term representations, such
as~\cite{Bodini2018,lmcs:848,10.1093/logcom/exx018} admit a common, general
symbolic specification template, differing only in weights assigned to specific
constructors (\emph{i.e.}~variables, applications, and abstractions). Let $\CS{L}$,
$\CS{F}$, $\CS{B}$, $\CS{A}$, and $\CS{U}$ denote the class of (open or closed)
$\lambda$\nobreakdash-terms, free and bound variables, a single application node, and
a single abstraction node, respectively. Then, $\CS{L}$ admits the following specification:
\begin{equation}\label{eq:random:generation:specification}
    \CS{L} = \CS{F} + (\CS{A} \times {\CS{L}}^2) + \CS{U} \times
    \text{subs}(\CS{F} \mapsto \CS{F} + \CS{D}, \CS{L}).
\end{equation}

In words, a $\lambda$\nobreakdash-term is either a free variable in $\CS{F}$, an
application of two terms $(\CS{A} \times {\CS{L}}^2)$ joined by an application
node, or an abstraction followed by a lambda term in which some free variables
become bound by the head abstraction (\emph{i.e.}~are substituted by some bound
variables in $\CS{D}$). Since not all free variables have to be bound by the
topmost abstraction, free variables can also be substituted for variables in
$\CS{F}$.

Let $FV(T)$ denote the set of free variables in $T$. Consider the following
\emph{bivariate} generating function $L(z, f)$ defined as
\begin{equation}
    L(z, f) = \sum_{T \in \CS{L}} z^{|T|} f^{|FV(T)|}.
  \end{equation}
$L(z,f)$ sums over all terms $T$ using variable $z$ to \emph{account} for $T$s
size, and variable $f$ to account for the number of its free variables. For
instance, consider the combinator $\mathbf{S} =\lambda x y z. x z (y z)$. If we
assume that $|\mathbf{S}| = 10$ (one size unit for each abstraction,
applications, and variable), then its respective monomial in $L(z,f)$ becomes
$z^{10}f^{0}$.

If we group matching monomials we note the coefficient $[z^nf^k]L(z, f)$ stands
for the number of terms of size $n$ and exactly $k$ free variables. In
particular, the coefficient $[z^n f^0]L(z,f)$ corresponds to the number of
closed $\lambda$\nobreakdash-terms of size $n$. We can also look at whole series
associated with specific powers. For instance, $[f^0]L(z,f) = \sum_{n \geq 0}
z^n f^0$ corresponds to the univariate generating function enumerating closed
$\lambda$\nobreakdash-terms. In general, $[f^k]L(z,f)$ is the generating
function corresponding to $\lambda$\nobreakdash-terms with exactly $k$ free
variables.

Unfortunately, by~\eqref{eq:random:generation:specification} the series
$[f^0]L(z,f)$ depends on $[f^1]L(z,f)$ which, in turn, depends on $[f^2]L(z,f)$,
\emph{etc.} Note that each new abstraction can bind occurrences of a
variable. Consequently, the abstraction body can have one more free variable than the
whole term. Due to this apparent infinite inductive nesting, finding the
generating function $[f^0]L(z,f)$ for closed terms, as well as generating
respective terms using Boltzmann models, is as difficult as the general case
$[f^k]L(z, f)$.

\begin{remark}
    Despite the inadmissible
    specification~\eqref{eq:random:generation:specification} it is still
    possible to use recursive samplers for closed $\lambda$\nobreakdash-terms.
    The number of free variables is bounded by the overall term size and, in
    particular, related to the number of term applications. Hence, for target
    term size $n$ we essentially need to memorise all of the values $[z^{\leq
    k}f^0]L(z,f),\ldots,[z^{\leq k},f^n]L(z,f)$ for $k = 0,\ldots,n$.

    Still, such a method allows us to sample terms of relatively small sizes,
    especially if the involved numbers $[z^{\leq k}f^m]L(z,f)$ grow too fast,
    cf.~\cite{Bodini2018,Wang05generatingrandom}. Consequently, the standard
    practice is to \emph{restrict} the class of generated terms, using
    practically justifiable criteria, such as limiting their unary height, number
    of abstractions, or the maximal allowed De~Bruijn index, see
    \emph{e.g.}~\cite{BODINI201845,Bodini2018,BendkowskiThesis}. In all these cases,
    the infinite inductive nesting of parameters becomes finite, and so it is
    possible to use more efficient sampling techniques, such as Boltzmann
    sampling.
\end{remark}

Under certain term representations, it is possible to sample large, unrestricted
closed $\lambda$\nobreakdash-terms, despite the infinite nesting problem. Let us
consider the class of $\lambda$\nobreakdash-terms in the De~Bruijn notation. For
convenience, assume the natural size notion in which each constructor (including
the successor and zero) weights one~\cite{10.1093/logcom/exx018}. Let $L(z)$
denote the corresponding generating function for plain (open or closed) terms and $L_m(z)$ denote
the generating function corresponding to so-called $m$\nobreakdash-open
$\lambda$\nobreakdash-terms, \emph{i.e.}~terms which when prepended with $m$ head
abstractions become closed. Note that if a term is $m$\nobreakdash-open, then it
is also $(m+1)$\nobreakdash-open.

With the help of $L_m(z)$ we can
rewrite~\eqref{eq:random:generation:specification} into a more verbose, infinite
specification exhibiting the precise relations among various openness levels:
\begin{align}\label{eq:random:generation:specification:2}
  \begin{split}
    L_0(z) &= z {L_0(z)}^2 + z L_1(z)\\
    L_1(z) &= z {L_1(z)}^2 + z L_2(z) + z\\
    L_2(z) &= z {L_2(z)}^2 + z L_3(z) + z + z^2\\
        \ldots\\
    L_m(z) &= z {L_m(z)}^2 + z L_{m+1}(z) + z \frac{1- z^m}{1-z}\\
        \ldots
    \end{split}
\end{align}
In words, if a term is $m$\nobreakdash-open, then it is either an application of
two $m$\nobreakdash-open terms, accounted for in the expression $z {L_m(z)}^2$,
an abstraction followed by an $(m+1)$\nobreakdash-term, accounted for in the
expression $z L_{m+1}(z)$ or, finally, one of $m$ available indices
$\idx{0},\ldots,\ldots, \idx{m-1}$. Note that $\idx{n} = S^{(n)} 0$ and so $|\idx{n}|
= n + 1$. Hence the final expression $z \frac{1- z^m}{1-z}$.

Although infinite, such a specification can be effectively analysed using
recently developed analytic
techniques~\cite{BODINI201845,DBLP:journals/corr/abs-1805-09419}. As a
by-product of this analysis, it is possible to devise a simple and highly
effective rejection-based sampling scheme for closed
$\lambda$\nobreakdash-terms.
Consider the following \emph{truncated} variant of~\eqref{eq:random:generation:specification:2}:
\begin{align}\label{eq:random:generation:specification:3}
    \begin{split}
        L_{0,N}(z) &= z {L_{0,N}(z)}^2 + z L_{1,N}(z)\\
        L_{1,N}(z) &= z {L_{1,N}(z)}^2 + z L_{2,N}(z) + z\\
        L_{2,N}(z) &= z {L_{2,N}(z)}^2 + z L_{3,N}(z) + z + z^2\\
        \ldots\\
        L_{N,N}(z) &= z {L_{N,N}(z)}^2 + z L(z) + z \frac{1- z^N}{1-z}\\
    L(z) &= z {L(z)}^2 + z L(z) + \frac{z}{1-z}
    \end{split}
\end{align}
Note that instead of unfolding the definition of $L_{m+1}(z)$ indefinitely, we
stop at level $N$ and use the generating function $L(z)$ corresponding to
\emph{all} $\lambda$\nobreakdash-terms at the final level $N$. Such a
finite specification defines a subclass of plain $\lambda$\nobreakdash-terms
excluding open terms without sufficiently large index values. In particular,
this class of terms approximates quite well the set of closed terms.

Following~\cite{BODINI201845} we
use~\eqref{eq:random:generation:specification:3} to sample plain
$\lambda$\nobreakdash-terms and reject those which are not closed. The
asymptotic number $[z^n]L_{0,N}(z)$ of such terms approaches the asymptotic
number $[z^n]L_0(z)$ of closed terms exponentially fast as $N$ tends to
infinity. In consequence, the imposed rejection does not influence the linear
$O(n)$ time complexity of the sampler drawing objects from some interval
$[(1-\varepsilon)n, (1+\varepsilon)n]$. In fact, with $N = 20$ the probability
of rejection is already of order $10^{-8}$, cf.~\cite[Section
5.4]{BODINI201845}. Consequently, such a sampler rarely needs to dismiss
generated terms and is able to generate large random closed
$\lambda$\nobreakdash-terms of sizes in the millions.

\begin{remark}
The discussed sampling method works for a broad class of size notions
    within the De~Bruijn term representation~\cite{GittenbergerGolebiewskiG16}.
    It is, however, specific to the unary encoding of De~Bruijn indices.

The intuitive reason behind the success of the above rejection sampler lies is
the fact that random plain terms strongly resemble closed ones. Typically, large
random terms are \emph{almost} closed. For instance, on average, a random term
contains a constant number of free variables, and needs just a few additional
head abstractions to become closed~\cite{DBLP:journals/corr/abs-1805-09419}.
\end{remark}

\medskip
{\bf Typeable $\lambda$-terms.}
Much like closed $\lambda$\nobreakdash-terms, simply-typed terms do not admit a
finite, admissible specification from which we could derive efficient samplers.
Simple recursive samplers work, however cease to be effective already for small
term sizes, cf.~\Cref{sec:random:generation:remark:lescanne:grygiel}. Sampling
simply-typed terms seems notoriously more challenging than sampling closed ones.
Even rejection sampling, whenever applicable, admits serious limitations due to
the imminent \emph{asymptotic sparsity problem} --- asymptotically almost no
term, be it either plain or closed, is at the same time (simply) typeable. This
problem is not only intrinsic to De~Bruijn models,
see~\cite{10.1007/978-3-662-49192-8_15}, but also seems to hold in models in
which variables contribute constant weight to the term size, see
\emph{e.g.}~\cite{lmcs:848,Bodini2018}.

Asymptotic sparsity of simply-typed $\lambda$\nobreakdash-terms is an
impenetrable barrier to rejection sampling techniques. As the term size tends to
infinity, so does the induced rejection overhead. In order to postpone this
inevitable obstacle, it is possible to use dedicated mechanisms interrupting the
sampler as soon as it is clear that the partially generated term cannot be
extended to a typeable one. The current state-of-the-art samplers take this
approach, combining Boltzmann models with modern logic programming execution
engines backed by highly-optimised unification
algorithms~\cite{bendkowski_grygiel_tarau_2018}. Nonetheless, even with these
sophisticated optimisations, such samplers are not likely to generate terms of
sizes larger than one hundred.

\begin{remark}
For some restricted classes of typeable terms there exist effective sampling
    methods, such as the already mentioned samplers for linear and affine
    terms~\cite{BODINI2013227}. Since Boltzmann models support all kinds of
    unambiguous context-free languages, it is also possible to sample other
    fragments of minimal intuitionistic logic, for instance, using finite
    truncations \emph{à la} Takahashi \emph{et al.}~\cite{TAKAHASHI1996144}.
\end{remark}

\subsection{Multiparametric samplers}\label{sec:multiparametric:samplers}
For some practical applications, especially in industrial property-based
software testing~\cite{7107466}, the uniform distribution of outcome structures
might not be the most useful choice. In fact, it can be argued that
most software bugs are miniscule, neglected corner cases, which will not be
caught using large, typical instances of random data,
see~\cite{Palka:2011:TOC:1982595.1982615,Runciman:2008:SLS:1411286.1411292}.

Remarkably, combinatorial samplers based on either Boltzmann models or the
recursive method, can be effectively \emph{tuned} so to distort the intrinsic
outcome distribution~\cite{doi:10.1137/1.9781611975062.9}. For instance,
consider again the class $\CS{L}$ of $\lambda$\nobreakdash-terms in the De~Bruijn
notation with natural size. The corresponding class $\CS{D}$ of De~Bruijn
indices satisfies
\begin{equation}
    \CS{D} = \CS{Z} \times \Seq{(\CS{Z})}.
\end{equation}
In words, a De~Bruijn index is a (possibly empty) sequence of
successors applied to zero. Both each successors and zero contribute weight one
to the overall term size. 

Such a specification inevitably dictates a geometric distribution of index
values and, consequently, an average constant distance between encoded variables
and their binders~\cite{DBLP:journals/corr/abs-1805-09419}. If such a
distribution is undesired, we can, for instance, force the sampler to output
terms according to a slightly more \emph{uniform distribution}. Suppose that we
change the De~Bruijn index specification to the more verbose
\begin{equation}
    \CS{D} = \CS{U}_0 \CS{Z} + \CS{U}_1 \CS{Z}^2 + \cdots + \CS{U}_k \CS{Z}^{k+1} +
    \CS{Z}^{k+2} \Seq{(\CS{Z})}
\end{equation}
assigning \emph{marking} variables $\CS{U}_0,\ldots,\CS{U}_k$ to the initial
$k+1$ distinct indices $\idx{0},\ldots,\idx{k}$, leaving the remaining indices
unmarked.  In doing so, we end up with a multiparametric specification for
$\lambda$\nobreakdash-terms where $\CS{Z}$, as usual, marks the term size, and
$\CS{U}_0,\ldots,\CS{U}_k$ mark some selected indices.

This new multivariate specification can be effectively \emph{tuned} in such a
way, that branching probabilities governing the sampler's decisions impose a
different, non-geometric distribution of index values. For that purpose the
problem is effectively expressed as a convex optimisation problem. With the
seminal advancements in interior-point methods for convex
programming~\cite{nesterov1994interior} the time required to compute the
branching probabilities is proportional to a polynomial involving the
optimisation problem size and the number of variables.

To illustrate the distortion effect, suppose that we mark the first nine indices
and impose a uniform distribution of $8\%$ of the total term size among them.
Remaining indices, are left
unaltered.~\Cref{fig:random:generation:lambda:term:freqs}, taken
from~\cite{DBLP:journals/corr/abs-1805-09419}, provides some empirical
frequencies obtained using a tuned sampler and its regular, undistorted
counterpart. Specific values were obtained for terms of size around $10,000$.

\begin{table*}[ht!]
\scalebox{.85}{
\begin{tabular}{c | c | c | c | c | c | c | c | c | c}
Index & $\idx{0}$ & $\idx{1}$ & $\idx{2}$ & $\idx{3}$ & $\idx{4}$ & $\idx{5}$ &
$\idx{6}$ & $\idx{7}$ & $\idx{8}$\\ \hline
Tuned frequency & $8.00\%$ & $8.00\%$ & $8.00\%$ & $8.00\%$ & $8.00\%$ &
$8.00\%$ & $8.00\%$ & $8.00\%$ & $8.00\%$\\
Observed frequency & $7.50\%$ & $7.77\%$ & $8.00\%$ & $8.23\%$ & $8.04\%$ &
$7.61\%$ & $8.53\%$ & $7.43\%$ & $9.08\%$\\
Default frequency & $21.91\%$ & $12.51\%$ & $5.68\%$ & $2.31\%$ &
$0.74\%$ & $0.17\%$ & $0.20\%$ & $0.07\%$ & - - -
\end{tabular}}
\caption{Empirical frequencies  of index distribution in relation to the term
    size.}
\label{fig:random:generation:lambda:term:freqs}
\end{table*}

\begin{remark}
    Multiparametric tuning requires a few mild technical assumptions, such as
the \emph{feasibility} of requested tuning (\emph{e.g.}~we cannot ask for the impossible).
It can also be used for exact-size sampling or even more involved admissible
constructions, including Pólya structures, labelled sets, \emph{etc}. The tuning
procedure can be carried out automatically and therefore used to rigorously
control the sample distribution without resorting to heuristics, or in-depth
expertise of the user.
\end{remark}

\begin{remark}
  Boltzmann samplers provide an incredibly versatile and efficient framework for
combinatorial generation. Despite the underlying use of heavy theories such as
analytic combinatorics and complex analysis, their \emph{design} and
\emph{usage} follow simple, recursive rules resembling a small, dedicated
\emph{domain-specific language}. With effective tools of their design and
tuning, Boltzmann samplers can be automatically compiled without additional
expertise in analytic combinatorics or complex
analysis~\cite{Duchon:2004:BSR:1024662.1024669,Canou:2009:FSR:1596627.1596637,doi:10.1137/1.9781611975062.9,PIVOTEAU20121711}.
We refer the more \emph{code-minded} reader to concrete
implementations\footnote{\url{https://github.com/Lysxia/boltzmann-samplers},\\
\url{https://github.com/maciej-bendkowski/paganini},\\
\url{https://github.com/maciej-bendkowski/boltzmann-brain},\\
\url{https://github.com/maciej-bendkowski/lambda-sampler}.}.

Let us note that analytic methods of (multivariate) Boltzmann sampling are not
the only rigorous methods of generating random algebraic data types or
$\lambda$\nobreakdash-terms in particular. While other methods exists, for
instance branching processes of Mista \emph{et al.}~\cite{mista}, analytic
methods seem to support the generation of far larger structures.
\end{remark}

\section{Conclusion}
\addtocontents{toc}{\protect\setcounter{tocdepth}{1}}
Random generation of $\lambda$\nobreakdash-terms is a young and active research
topic on the border of software testing, theoretical computer science and
combinatorics. It combines utmost practical experimentation and theoretical
advances in analytic combinatorics. Still, many important problems remain open.
For instance, how to deal with representations in which corresponding generating
functions cease to be analytic, or how to generate large random simply-typed
terms?

\subsection*{Acknowledgements}
We would like to thank Pierre Lescanne and Sergey Dovgal for encouraging us to
write this survey and their valuable remarks during its formation.

\printbibliography
\end{document}